\documentclass[preprint]{elsarticle}

\usepackage{lineno,hyperref}
\usepackage{hyperref}
\usepackage{amssymb,amsmath}

\journal{arXiv}

\newtheorem{thm}{Theorem}%[section]
\newtheorem{prop}{Proposition}%[section]
\newproof{pf}{Proof}

\renewcommand{\leq}{\leqslant}
\renewcommand{\geq}{\geqslant}

\def\mc{\mathcal}
\def\mb{\mathbb}
\def\ve{\varepsilon}

\def\ce{c}
\def\ra{\rightarrow}

\providecommand{\abs}[1]{\lvert#1\rvert}

\DeclareMathOperator{\Span}{Span}
\DeclareMathOperator{\Ker}{Ker}

\DeclareMathOperator{\im}{Im}

\DeclareMathOperator{\id}{Id}

\begin{document}

\begin{frontmatter}

\title{On quasi-Clifford Osserman curvature tensors}

%\author{Vladica Andreji\'c\fnref{myfootnote}}

\author[matf]{V.~Andreji\'c}
\ead{andrew@matf.bg.ac.rs}

\author[matf]{K.~Luki\'c}
\ead{katarina\_lukic@matf.bg.ac.rs}

\address[matf]{Faculty of Mathematics, University of Belgrade, Belgrade, Serbia}




\begin{abstract}
We consider pseudo-Riemannian generalizations of Osserman, Clifford, and the duality principle properties
for algebraic curvature tensors and investigate relations between them.
We introduce quasi-Clifford curvature tensors using a generalized Clifford family and show that they are Osserman.
This allows us to discover an Osserman curvature tensor that does not satisfy the duality principle.
We give some necessary and some sufficient conditions for the total duality principle.
\end{abstract}

\begin{keyword}
Osserman manifold\sep Clifford family\sep duality principle
\MSC[2010] Primary 53B30; Secondary 53C50
\end{keyword}

\end{frontmatter}

%\linenumbers

\section{Introduction}

Let $(\mc{V},g)$ be a (possibly indefinite) scalar product space of dimension $n$.
The squared norm of a vector $X\in\mc{V}$ is the real number $\ve_X=g(X,X)$.
The sign of the squared norm distinguishes all vectors $X\in\mc{V}$ into three different types.
A vector $X\in\mc{V}$ is spacelike if $\ve_X>0$; timelike if $\ve_X<0$; null if $\ve_X=0$.
Especially, a vector $X\in \mc{V}$ is nonnull if $\ve_X\neq 0$ and it is unit if $\ve_X\in \{-1,1\}$.

An algebraic curvature tensor on $(\mc{V},g)$ is a quadri-linear map $R\colon \mc{V}^4\ra\mb{R}$ that 
satisfies usual $\mathbb{Z}_2$-symmetries and the first Bianchi identity.
In the presence of an orthonormal basis $(E_1,\dots,E_n)$ in $\mc{V}$, 
we have the associated Jacobi operator $\mc{J}_X\colon \mc{V}\ra\mc{V}$ for $X\in\mc{V}$ by 
$\mc{J}_X(Y)=\sum_{i=1}^n \ve_{E_i} R(Y,X,X,E_i)E_i$.
The Jacobi operator is a self-adjoint endomorphism on $\mc{V}$, and therefore it is diagonalizable if $g$ is definite.
However, this is no longer true in the indefinite setting, so if $\mc{J}_X$ is diagonalizable for any nonnull $X$ we say that $R$ is Jacobi-diagonalizable.
In general, the eigen-structure of $\mc{J}_X$ is determined by the Jordan normal form (the number and the sizes of the Jordan blocks).

We say that $R$ is timelike Osserman (or spacelike Osserman) if the characteristic polynomial of the Jacobi operator $\mc{J}_X$
is independent of unit timelike (or spacelike) $X\in\mc{V}$.
We say that $R$ is timelike Jordan-Osserman (or spacelike Jordan-Osserman) if the Jordan normal form of $\mc{J}_X$ 
is independent of unit timelike (or spacelike) $X$.
An algebraic curvature tensor is Osserman if it is both timelike and spacelike Osserman,
and it is Jordan-Osserman if it is both timelike and spacelike Jordan-Osserman.
It is known that spacelike Osserman and timelike Osserman are equivalent properties, and $R$ is Osserman
if and only if $\det(\ve_X\lambda\id-\mc{J}_X)=0$ is the same equation for all nonnull $X$.
However, spacelike Jordan-Osserman and timelike Jordan-Osserman are not equivalent (see Section \ref{jd}).

In the Riemannian setting ($g$ is positive definite), one of important features of an Osserman 
algebraic curvature tensor is the duality principle, given by Raki\'c \cite{Ra2}.
Generalizations to a pseudo-Riemannian setting (see Andreji\'c and Raki\'c \cite{AR1,AR2}) 
are possible by the following implication,
\begin{equation}\label{dual}
Y \text{ is an eigenvector of } \mc{J}_X\quad \Longrightarrow \quad X \text{ is an eigenvector of } \mc{J}_Y.
\end{equation}
We shall use two kinds of duality depending on our $(X,Y)$ domain.
If \eqref{dual} holds for all $X,Y\in\mc{V}$ with the restriction $\ve_X\neq 0$, we say that $R$ is Jacobi-dual ($R$ satisfies the duality principle),
and if the only restriction is $X\neq 0$, we say that $R$ is totally Jacobi-dual (see Andreji\'c and Raki\'c \cite{AR2}).

After Raki\'c established the duality principle and proved that a Riemannian Osserman algebraic curvature tensor is Jacobi-dual (see \cite{Ra2}),
it is extensively studied in the pseudo-Riemannian settings. 
The best results were recently given by Nikolayevsky and Raki\'c \cite{NR2}, where they showed that any Jordan-Osserman $R$ is Jacobi-dual.
They also proved that if the set of those $X\in\mc{V}$ for which $\mc{J}_X$ is diagonalizable has a nonempty interior ($R$ is semisimple)
then $R$ is Osserman if and only if it is Jacobi-dual.

Additionally, we have the following partial results.
Any four-dimensional Osserman $R$ is Jacobi-dual, see Andreji\'c \cite{A3}.
Any Lorentzian totally Jacobi-dual $R$ has constant sectional curvature, see Andreji\'c and Raki\'c \cite{AR2}.
Any four-dimensional Jacobi-dual $R$ such that $\mc{J}_X$ is diagonalizable for some nonnull $X$ is Osserman, see Andreji\'c \cite{A7}.

\section{Quasi-Clifford curvature tensors} \label{qcct}

A tensor of constant sectional curvature $1$ is a very first example of an algebraic curvature tensor,
\begin{equation*}
R^0(X,Y,Z,W)=g(Y,Z)g(X,W)-g(X,Z)g(Y,W).
%\mc{R}^0(X,Y)Z =g(Y,Z)X-g(X,Z)Y.
\end{equation*}
Additionally, any skew-adjoint endomorphism $J$ on $\mc{V}$ generates a new example by 
\begin{equation*}
R^J(X,Y,Z,W) = g(JX,Z)g(JY,W) -g(JY,Z)g(JX,W) +2g(JX,Y)g(JZ,W).
\end{equation*}
Therefore, a linear combination 
\begin{equation}\label{clifford}
R=\mu_0 R^0+\sum_{i=1}^m \mu_i R^{J_i}
\end{equation}
is an algebraic curvature tensor for skew-adjoint endomorphisms $J_1,\dots,J_m$ on $\mc{V}$.

A Clifford family is an anti-commutative family of skew-adjoint complex structures $J_i$, $1\leq i\leq m$.
Algebraic curvature tensors of form \eqref{clifford} associated with a Clifford family were introduced by Gilkey \cite{Gi3,GSV}.
However, it is natural to consider the generalization, an anti-commutative family of $J_i$ such that $J_i^2=\ce_i\id$
for some $\ce_i\in\mb{R}$, that is, the Hurwitz-like relations,
\begin{equation}\label{hurwitz}
J_i J_j+J_j J_i=2\ce_i\delta_{ij}\id,
\end{equation}
hold for $1\leq i,j\leq m$.
We say that an algebraic curvature tensor $R$ is quasi-Clifford if it has a form \eqref{clifford} with the Hurwitz-like relations \eqref{hurwitz}.
Especially, $R$ is Clifford if it is quasi-Clifford with $\ce_i=-1$ for all $1\leq i\leq m$.

It is well known that a Clifford algebraic curvature tensor is Osserman.
However, according to Nikolayevsky \cite[Section 2]{N0}, in the Riemannian setting the converse is true (an Osserman $R$ is Clifford)
in all dimensions except $n=16$, and also in many cases when $n=16$.
The only known (Riemannian) counterexample is the curvature tensor $R^{OP^2}$ of the Cayley projective plane,
or more precisely, any algebraic curvature tensor of the form $\mu R^{OP^2}+\xi R^1$, where $R^1$ is the curvature tensor
of the unit sphere $S^{16}(1)$ and $\mu\neq 0$.

Let us start with a quasi-Clifford $R$ and an arbitrary vector $X\in\mc{V}$.
Each $J_i$ is skew-adjoint which implies $g(J_iX,X)=0$ and simplifies the calculation of the Jacobi operator, 
\begin{equation}\label{jacobi}
\mc{J}_X(Y)=  \mu_0 (\ve_X Y-g(Y,X)X) - 3\sum_{i=1}^m \mu_i g(Y,J_iX)J_iX.
\end{equation}
Additionally, the equation \eqref{hurwitz} implies $g(J_i X,J_j X)=0$ and $\ve_{J_iX}=-\ce_i\ve_X$ for $1\leq i\neq j\leq n$.
If we denote $\mc{F}_t= \{X,J_1 X,\dots,J_t X\}$ for $1\leq t\leq m$, then we obtain
\begin{equation*}
\mc{J}_X(J_iX)=\ve_X(\mu_0+3\ce_i\mu_i)J_iX,\quad \mc{J}_X(Z)=\ve_X\mu_0 Z,
\end{equation*}
for all $1\leq i\leq m$ and $Z\in\mc{F}_m^{\perp}$.
It is important to distinguish the case $\ce_i\neq 0$, $1\leq i\leq k$ from the case $\ce_i=0$, $k<i\leq m$.

For a nonnull $X$, the set $\mc{F}_k$ consists of mutually orthogonal nonnull eigenvectors.
Thus, $\Span\mc{F}_k$ and $\mc{F}_k^{\perp}$ are nondegenerate, and we can consider the restriction $\widetilde{\mc{J}}_X$ of $\mc{J}_X$ to $\mc{F}_k^{\perp}$
as a self-adjoint endomorphism on $\mc{F}_k^{\perp}$.
Eigenvectors corresponding to distinct eigenvalues of a self-adjoint endomorphism are mutually orthogonal,
including complex eigenvectors for complex eigenvalues from the complexification $\mc{V}^{\mb{C}}\cong \mc{V}\oplus i\mc{V}$.
However, all vectors orthogonal to $\mc{F}_m\setminus \mc{F}_k$ are eigenvectors with the eigenvalue $\ve_X\mu_0$,
and consequently $\widetilde{\mc{J}}_X$ has no other eigenvalues.
Hence,
\begin{equation*}
\det(\ve_X\lambda\id-\mc{J}_X)=  (\ve_X)^n\lambda(\lambda-\mu_0)^{n-k-1}\Pi_{i=1}^k(\lambda-(\mu_0+3\ce_i\mu_i)),
\end{equation*}
which means that $R$ is Osserman.

\begin{thm}
Any quasi-Clifford algebraic curvature tensor is Osserman.
\end{thm}

Additionally, the Jordan normal form of $\mc{J}_X$ has the critical part on $\mc{F}_k^{\perp}$, where for $\widetilde{\mc{J}}_X$ we have 
\begin{equation*}
\im(\widetilde{\mc{J}}_X-\ve_X\mu_0\id)\subseteq \Span\{J_{k+1}X,\dots,J_mX\}\subseteq \Ker(\widetilde{\mc{J}}_X-\ve_X\mu_0\id),
\end{equation*}
and therefore $\widetilde{\mc{J}}_X-\ve_X\mu_0\id$ is two-step nilpotent, $(\widetilde{\mc{J}}_X-\ve_X\mu_0\id)^2=0$.
Thus, a quasi-Clifford $R$ allows only Jordan blocks of size $2$.

However, let us remark that a pseudo-Riemannian manifold $(\mb{R}^4,g)$ induced by the metric 
$g=x_2x_3 dx_1^2 - x_1x_4 dx_2^2+ 2 dx_1dx_2 + 2 dx_1dx_3 + 2 dx_2dx_4$,
at any point has the curvature tensor that is Jordan-Osserman such that the Jordan normal form has a Jordan block of size $3$ (see \cite[Remark 4.1.2]{GKL}). 
This means that the converse question fails in the signature $(2,2)$, where an Osserman $R$ is not necessarily quasi-Clifford. 

\section{Jacobi-duality}\label{jd}

We shall follow and generalize the arguments from Andreji\'c and Raki\'c \cite{AR2} with the purpose 
to investigate whether a quasi-Clifford $R$ is Jacobi-dual.
Let $X\in\mc{V}$ be nonnull and suppose that $Y$ is an eigenvector of $\mc{J}_X$, that is, $\mc{J}_X(Y)=\ve_X\lambda Y$ holds for some $\lambda\in\mb{R}$.
Then \eqref{jacobi} implies
\begin{equation}\label{ipsilon}
\ve_X(\lambda-\mu_0)Y=  -\mu_0  g(Y,X)X - 3\sum_{i=1}^m \mu_i g(Y,J_iX)J_iX,
\end{equation}
while by interchanging the roles of $X$ and $Y$ in \eqref{jacobi} we have
\begin{equation}\label{jacobi2}
\mc{J}_Y(X)=  \mu_0 (\ve_Y X-g(X,Y)Y) - 3\sum_{i=1}^m \mu_i g(X,J_iY)J_iY.
\end{equation}
If $\lambda\neq\mu_0$, then we can express $Y$ from \eqref{ipsilon} and get
%\begin{equation}\label{ipsilon1}
%Y= \frac{-\mu_0 g(Y,X)}{\ve_X (\lambda-\mu_0)} X - 3\sum_{i=1}^m \frac{\mu_i g(Y,J_iX)}{\ve_X (\lambda-\mu_0)} J_iX,
%\end{equation}
%\begin{align*}
%\mc{J}_Y(X)=&  \mu_0 \left(\ve_Y X-g(X,Y)\left(\frac{-\mu_0 g(Y,X)}{\ve_X (\lambda-\mu_0)} X - 
%3\sum_{j=1}^m \frac{\mu_j g(Y,J_jX)}{\ve_X (\lambda-\mu_0)} J_jX\right)\right) \\
%&- 3\sum_{i=1}^m \mu_i g(X,J_iY)J_i\left(\frac{-\mu_0 g(Y,X)}{\ve_X (\lambda-\mu_0)} X - 
%3\sum_{j=1}^m \frac{\mu_j g(Y,J_jX)}{\ve_X (\lambda-\mu_0)} J_jX\right).
%\end{align*}
%It gives
\begin{align*}
\mc{J}_Y(X)=&\phantom{+}  \mu_0 \left(\ve_Y + \frac{\mu_0 (g(X,Y))^2}{\ve_X (\lambda-\mu_0)}\right) X\\
&+ \frac{3\hspace{.2mm}\mu_0g(X,Y)}{\ve_X (\lambda-\mu_0)}\sum_{i=1}^m \mu_i \left(g(Y,J_iX)+g(X,J_iY)\right) J_iX  \\
&+ \frac{9}{\ve_X (\lambda-\mu_0)}\sum_{i=1}^m \sum_{j=1}^m \mu_i\mu_j g(X,J_iY)g(Y,J_jX) J_iJ_jX\\
=&\left(\mu_0 \ve_Y + \frac{\mu_0^2 (g(X,Y))^2}{\ve_X (\lambda-\mu_0)} - 
\frac{9}{\ve_X (\lambda-\mu_0)} \sum_{i=1}^m \mu_i^2 (g(Y,J_iX))^2 \ce_i\right) X,
\end{align*}
so $X$ is an eigenvector of $\mc{J}_Y$.
Otherwise, $\lambda=\mu_0$ gives
\begin{equation}\label{lin}
\mu_0  g(Y,X)X + 3\sum_{i=1}^m \mu_i g(Y,J_iX)J_iX=0.
\end{equation}
If the set $\mc{F}_m$ is linearly independent, we have $g(Y,X)=g(Y,J_iX)=g(J_iY,X)=0$,
and therefore $\mc{J}_Y(X)=  \ve_Y \mu_0 X$, and again, $X$ is an eigenvector of $\mc{J}_Y$.
For a nonnull $X$, $\Span\mc{F}_k$ is a nondegenerate subspace of $\mc{V}$ orthogonal to $\mc{F}_m\setminus \mc{F}_k$.
Thus, if $m=k$ or $m=k+1$, then the set $\mc{F}_m$ is linearly independent.

\begin{thm}
Any quasi-Clifford algebraic curvature tensor with at most one $c_i=0$ is Jacobi-dual.
\end{thm}
 
However, this is no longer true if there are at least two $J_i$ with $c_i=0$.
Let $(T_1,\dots T_p,S_1,\dots,S_q)$ be an orthonormal basis in a scalar product space $(\mc{V},g)$ of signature $(p,q)$.
Let us set an endomorphism $J$ on $\mc{V}$ by
\begin{align*}
JT_1=T_2+S_2=-JS_1,\quad -JT_2=T_1+S_1=JS_2\\
JT_3=T_4+S_4=-JS_3,\quad -JT_4=T_3+S_3=JS_4\\
JT_5=\dots=JT_p=JS_5=JS_6=\dots JS_q=0.
\end{align*}
It is easy to check that $J$ is skew-adjoint with $J^2=0$.
In the case $4=p<q$, an Osserman algebraic curvature tensor $R=R^J$ is timelike Jordan-Osserman,
but it is not spacelike Jordan-Osserman, which is similar to Gilkey and Ivanova \cite{GI0} and  Gilkey \cite[Section 3.2]{Gi1}.
Let us introduce an additional endomorphism $K$ on $\mc{V}$ by
\begin{align*}
KT_1=T_2+S_2=-KS_1,\quad -KT_2=T_1+S_1=KS_2\\
KT_3=T_4+S_5=-KS_3,\quad -KT_4=T_3+S_3=KS_5\\
KT_5=\dots=KT_p=KS_4=KS_6=\dots KS_q=0.
\end{align*}
It is just changing roles of $S_4$ and $S_5$ in $J$, so $K$ is also skew-adjoint with $K^2=0$.
If we set $R=R^J-R^K$, then from \eqref{jacobi}, $\mc{J}_X(Y)=  3(g(Y,KX)KX-g(Y,JX)JX)$.
Since $JT_1=KT_1=T_2+S_2$ we have $\mc{J}_{T_1}(Y)=0$ for any $Y\in\mc{V}$.
Additionally, for $Y=T_2+\sqrt{2}S_4$ we have $g(Y,JT_1)=g(Y,KT_1)=g(T_2+\sqrt{2}S_4,T_2+S_2)=-1$, and therefore
$\mc{J}_Y(T_1)=  3(g(Y,JT_1)JY-g(Y,KT_1)KY)=-3\sqrt{2}(T_3+S_3)$. Thus, 
\begin{equation*}
\mc{J}_{T_1}(T_2+\sqrt{2}S_4)=0, \quad \mc{J}_{T_2+\sqrt{2}S_4}(T_1)=-3\sqrt{2}(T_3+S_3)
\end{equation*}
show that $R$ is not Jacobi-dual. 
Moreover, our counterexample contains mutually orthogonal unit vectors $X=T_1$ and $Y=T_2+\sqrt{2}S_4$,
such that \eqref{dual} does not work.
In this way we were able to discover an Osserman $R$ that is not Jacobi-dual.

\begin{thm}
There exist quasi-Clifford (and therefore Osserman) algebraic curvature tensors which are not Jacobi-dual.
\end{thm}

\section{Total Jacobi-duality}

Skew-adjoint endomorphisms with $J_i^2=0$ change the Jordan normal form of $\mc{J}_X$ and therefore they are inadequate 
for the duality principle, which we have already seen in the previous section.
Hence, we shall exclude them ($m=k$), which leaves only the $J_i$ that are automorphisms.
Without loss of generality, using the rescaled $(1/\sqrt{\abs{\ce_i}})J_i$, we can suppose $\ce_i\in\{-1,1\}$, and such $R$ we called semi-Clifford.
From Section \ref{qcct}, it is easy to see that semi-Clifford $R$ is Jacobi-diagonalizable and consequently Jordan-Osserman.

A semi-Clifford algebraic curvature tensor is generated by 
a family of anti-commutative skew-symmetric orthogonal and anti-orthogonal operators on $\mc{V}$.
In fact, these are complex structures ($\ce_i=-1$) and product structures ($\ce_i=1$).
It worth noting that a product structures $J_i$ change the signature because of $\ve_{J_iX}=-\ve_X$, 
so in a non-neutral signature any semi-Clifford $R$ is Clifford. 

We have already seen that a semi-Clifford $R$ is Jacobi-dual, and the next step is to investigate whether $R$ is totally Jacobi-dual.
The previous paper \cite{AR2} gives only a sufficient condition that $\mc{F}_m$ is linearly independent.
Namely, if $X$ is null, then the equation \eqref{jacobi} for $\mc{J}_X(Y)=0$ yields the equation \eqref{lin},
where the linear independence of $\mc{F}_m$, as before, gives $\mc{J}_Y(X)=  \ve_Y \mu_0 X$.
Therefore, we should examine whether the set $\mc{F}_m=\{X,J_1 X,\dots,J_m X\}$ is linearly independent for a null vector $X\neq 0$.

Let us suppose that
\begin{equation}\label{start}
\theta_0 X+\theta_1 J_1X+\dots+\theta_m J_mX=0
\end{equation}
holds for some $\theta_0,\theta_1,\dots,\theta_m\in\mb{R}$ and $X\neq 0$.
Applying an endomorphism 
\begin{equation*}
(-1)^{\alpha}J^\alpha=(-1)^{\alpha_1+\dots+\alpha_m}J_m^{\alpha_m}\dots J_1^{\alpha_1}
\end{equation*}
for $\alpha=(\alpha_1,\dots,\alpha_m)\in\{0,1\}^m$
we get a new equation 
\begin{equation*}
\sum_{i=0}^m (-1)^{\alpha}\theta_i J^\alpha J^{e_i}X=0,
\end{equation*}
%$\theta_0 J^\alpha X+\theta_1 J^\alpha J^{e_1}X+\dots+\theta_m J^\alpha J^{e_m}X=0$,
where $e_i=(\delta_{i1},\delta_{i2},\dots,\delta_{im})$ with additional $e_0=(0,\dots,0)$, i.e. $J^{e_i}=J_i$, and $J^{e_0}=\id$. 
It is easy to see that
\begin{equation*}
(-1)^{\alpha}J^\alpha J^{e_i}=(-1)^{\alpha_{i}+\dots+\alpha_m}(\ce_i)^{\alpha_i}J^{\alpha\pm e_i},
\end{equation*}
where ${\alpha\pm e_i}$ and $\alpha$ differ only in the $i$-th slot. 
Thus we have a homogeneous system of $2^m$ linear equations with $2^m$ unknowns,
\begin{equation}\label{system}
\sum_{\beta}M_{\alpha\beta}J^{\beta}X=0,
\end{equation}
with $M_{\alpha\alpha}=(-1)^{\alpha_{i}+\dots+\alpha_m}\theta_0$, 
$M_{\alpha(\alpha\pm e_i)}=(-1)^{\alpha_{i}+\dots+\alpha_m}(\ce_i)^{\alpha_i}\theta_i$ for $1\leq i\leq m$,
and $M_{\alpha\beta}=0$ otherwise.
Consider the matrix $M^2$, and calculate its entries,
\begin{align*}
(M^2)_{\alpha\alpha}
&=M_{\alpha\alpha}M_{\alpha\alpha}+ M_{\alpha(\alpha\pm e_1)}M_{(\alpha\pm e_1)\alpha}+\dots+M_{\alpha(\alpha\pm e_m)}M_{(\alpha\pm e_m)\alpha}\\
&=\theta_0^2-\ce_1\theta_1^2-\dots-\ce_m\theta_m^2,\\
(M^2)_{\alpha(\alpha\pm e_i)}
&=M_{\alpha\alpha}M_{\alpha(\alpha\pm e_i)}+M_{\alpha(\alpha\pm e_i)}M_{(\alpha\pm e_i)(\alpha\pm e_i)}=0,\\
(M^2)_{\alpha(\alpha\pm e_i\pm e_j)}
&=M_{\alpha(\alpha\pm e_i)}M_{(\alpha\pm e_i)(\alpha\pm e_i\pm e_j)}+M_{\alpha(\alpha\pm e_j)}M_{(\alpha\pm e_j)(\alpha\pm e_i\pm e_j)}=0,
\end{align*}
and $(M^2)_{\alpha\beta}=0$ otherwise. 
Thus, $M^2$ is a diagonal matrix with 
\begin{equation*}
(\det M)^2=\det(M^2)=(\theta_0^2-\ce_1\theta_1^2-\dots-\ce_m\theta_m^2)^{2^m}.
\end{equation*}
It is important to notice that if $\det M\neq 0$, the system \eqref{system} has the unique zero solution $X=J_1X=\dots=J_mX=0$, 
which is absurd since $X\neq 0$. Hence, we have the following statement.

\begin{thm}\label{tete}
If $\theta_0 X+\theta_1 J_1X+\dots+\theta_m J_mX=0$ holds for $X\neq 0$ then
\begin{equation}\label{det1}
\theta_0^2-\ce_1\theta_1^2-\dots-\ce_m\theta_m^2=0.
\end{equation}
\end{thm}

In the case of Clifford $R$, we have $\ce_i=-1$ for all $1\leq i\leq m$, so \eqref{det1} gives $\theta_0^2+\theta_1^2+\dots+\theta_m^2=0$,
and therefore $\theta_i=0$ holds for all $0\leq i\leq m$, $\mc{F}_m$ is linearly independent and consequently $R$ is totally Jacobi-dual.

\begin{thm}
Any Clifford algebraic curvature tensor is totally Jacobi-dual.
\end{thm}

However, there are some problems in the case that we have $\ce_i=1$ for some $i$.
If $R$ is semi-Clifford that is not totally Jacobi-dual, then there exists a pair of vectors $(X,Y)$
such that $X\neq 0$ is null with $\mc{J}_X(Y)=0$, where $\mc{J}_Y(X)$ is not proportional to $X$.
As before, we need \eqref{start}, $\theta_0 X+\theta_1 J_1X+\dots+\theta_m J_mX=0$,
such that $\theta_0=\mu_0 g(Y,X)$, and $\theta_i=3\mu_i g(Y,J_iX)$, $1\leq i\leq m$.
Thus,
\begin{equation*}
\theta_0^2=\mu_0 g(Y,\theta_0 X)=-\mu_0\sum_{i=1}^m \theta_i g(Y,J_iX)= -\mu_0\sum_{i=1}^m \frac{\theta_i^2}{3\mu_i}.
\end{equation*}
If we include Theorem \ref{tete}, the necessary conditions become
\begin{equation}\label{condition}
\theta_0^2= \sum_{i=1}^m \ce_i\theta_i^2= -\mu_0\sum_{i=1}^m \frac{\theta_i^2}{3\mu_i}.
\end{equation}
Hence 
\begin{equation*}
\sum_{i=1}^m(\ce_i+\frac{\mu_0}{3\mu_i})\theta_i^2=0,
\end{equation*} 
which implies the following theorem.

\begin{thm}\label{semic}
If $(3\ce_i\mu_i+\mu_0)\mu_i>0$ for all $1\leq i\leq m$ or $(3\ce_i\mu_i+\mu_0)\mu_i<0$ for all $1\leq i\leq m$,
then the associated semi-Clifford $R$ is totally Jacobi-dual.
\end{thm}

\section{Anti-Clifford curvature tensor}

In the case that $\ce_i=1$ holds for all $1\leq i\leq m$ we say that $R$ is anti-Clifford.
Then, Theorem \ref{tete} for the hypothetical $\theta_0=0$ gives 
$\theta_1^2+\dots+\theta_m^2=0$ which shows that the set $\{J_1X,\dots,J_mX\}$ is linearly independent, %for $J_1^2=\dots J_m^2=\id$, 
but $\mc{F}_m$ can be linearly dependent.

We say that a subspace $\mc{W}$ of an indefinite scalar product space $(\mc{V},g)$ is totally isotropic if it consists only of null vectors,
which implies that any two vectors from $\mc{W}$ are mutually orthogonal.
We need the next well known statement about an isotropic supplement of $\mc{W}$.% (for example, see Clark \cite[Theorem 6.2]{Cl}).

\begin{prop}\label{recip}
Let $\mc{W}\subset \mc{V}$ be a totally isotropic subspace with basis $N_1,\dots,N_k$.
Then there exist a totally isotropic subspace $\mc{U}$, disjoint from $\mc{W}$, with basis $M_1,\dots,M_k$,
such that $g(N_i,M_j)=\delta_{ij}$ holds for $1\leq i,j\leq k$.
\end{prop}

\begin{pf} 
The proof is by induction on $k$, where the case $k=0$ is trivial.
Let us set $\mc{P}=\Span\{N_1,\dots,N_{k-1}\}$.
Since $\Span\{N_k\}$ is not a subspace of $\mc{P}$, $\mc{P}^{\perp}$ is not a subspace of $\Span\{N_k\}^{\perp}$,
and there exists $X\in\mc{P}^{\perp}$ such that $g(X,N_k)\neq 0$.
Then 
\begin{equation*}
M_k=\frac{-\ve_X}{2(g(X,N_k))^2}N_k+\frac{1}{g(X,N_k)}X
\end{equation*} 
is null with $g(N_k,M_k)=1$. 
The subspace $\Span\{N_k,M_k\}$ is nondegenerate since it has an orthonormal basis $(N_k+M_k)/{\sqrt{2}},(N_k-M_k)/{\sqrt{2}}$.
By construction $\mc{P}$ is a subspace of the nondegenerate $\Span\{N_k,M_k\}^{\perp}$
and we apply the induction hypothesis to get $M_1,\dots,M_{k-1}$ with desired properties. \qed
\end{pf}

Since $J_1X,\dots,J_mX$ form a basis of a totally isotropic subspace of $\mc{V}$, 
by the previous proposition there exists a basis $\{M_1,\dots,M_m\}$ of an isotropic supplement, such that 
$g(J_iX,M_j)=\delta_{ij}$ and $g(M_i,M_j)=0$ hold for $1\leq i,j\leq m$.
Then 
\begin{equation*}
Z=\sum_{i=1}^m \frac{\theta_i}{3\mu_i}M_i
\end{equation*}
has the properties $\theta_i=3\mu_i g(Z,J_iX)$, and consequently, by \eqref{condition}, $\theta_0=\mu_0 g(Z,X)$.
Moreover, $Z+W$ has the same properties for any $W\in\{J_1X,\dots,J_mX\}^{\perp}$.

From \eqref{jacobi2}, we need such $Y$ that $-\theta_0 Y + \sum_{i=1}^m \theta_i J_iY$ is not proportional to $X$.
Therefore, we search for $Y=Z+W$ such that $\mc{K}(Z+W)$ is not proportional to $X$ where $\mc{K}=-\theta_0\id+\sum_{i=1}^m\theta_iJ_i$.
For any nonnull $D$, the set $\{D,J_1D,\dots,J_mD\}$ is linearly independent (mutually orthogonal nonnull vectors),
and therefore $\mc{K}(D)=0$ implies $\theta_0=\dots=\theta_m=0$ which is impossible. 
Thus $\mc{K}(D)\neq 0$ holds for all nonnull $D$.

Let us suppose that $\dim\mc{V}=n>2m$, which enables a nonnull vector $H$ from $\{J_1X,\dots,J_mX,M_1,\dots,M_m\}^{\perp}$.
We already have $\mc{K}(X)=-2\theta_0 X$.
If $R$ is totally Jacobi-dual, we have additional $\mc{K}(Z)=\zeta X$ and $\mc{K}(Z+H)=\xi X$.
Then $\mc{K}(H)=(\xi-\zeta) X\neq 0$, so $\mc{K}(Z-\frac{\zeta}{\xi-\zeta}H)=0$, which implies that $Z-\frac{\zeta}{\xi-\zeta}H$ is null,
and therefore $\zeta=0$.
Similarly, $\mc{K}(X+\frac{2\theta_0}{\xi}H)=0$, which implies that $X+\frac{2\theta_0}{\xi}H$ is null, and therefore $\theta_0=0$,
which is a contradiction.
Hence, we have the following theorem.

\begin{thm}\label{nontot}
If there exist $\theta_0,\dots,\theta_m\in\mb{R}$ (not all equal to zero), such that $\theta_0 X+\theta_1 J_1X+\dots+\theta_m J_mX=0$ holds for some nonnull $X$, 
with the condition \eqref{condition}, then the associated anti-Clifford algebraic curvature tensor of dimension $n>2m$ is not totally Jacobi-dual.
\end{thm}

In the end, let us show some examples of anti-Clifford algebraic curvature tensors which are not totally Jacobi-dual.

In the case $m=1$, Theorem \ref{semic} gives a necessary condition $\mu_0+3\mu_1=0$.
A skew-adjoint product structure $J$ given by $JT_i=S_i$ and $JS_i=T_i$ for all $1\leq i\leq t$, $n=2t\geq 4$,
where $(T_1,\dots,T_t,S_1,\dots,S_t)$ is an orthonormal basis in a scalar product space $(\mc{V},g)$ of neutral signature,
provides an anti-Clifford $R=3\mu R^0-\mu R^J$ for $\mu\neq 0$.
We can take $X=S_1+T_1$, because of the linear dependence $X=JX$, and apply Theorem \ref{nontot}.
%We can take $X=S_1+T_1$ and $Y=S_1+S_2$ to get
%\begin{equation*}
%\mc{J}_{S_1+T_1}(S_1+S_2)=0,\quad \mc{J}_{S_1+S_2}(S_1+T_1)=(\sqrt{2}-1)\mu_0(S_1+T_1)-\mu_0(S_2+T_2).
%\end{equation*}

In the case $m=2$, Theorem \ref{semic} gives $(\mu_0+3\mu_1)(\mu_0+3\mu_2)\mu_1\mu_2\leq 0$ as a necessary condition.
Consider skew-adjoint product structures $J$ and $K$ given by 
$JT_{2i-1}=S_{2i-1}$, $JT_{2i}=S_{2i}$, $KT_{2i-1}=S_{2i}$, and $JT_{2i}=-S_{2i-1}$, for all $1\leq i\leq t$, $n=4t\geq 8$,
where $(T_1,\dots,T_{2t},S_1,\dots,S_{2t})$ is an orthonormal basis in $(\mc{V},g)$.
We can choose $X=\cos\alpha\,T_1+\sin\alpha\,T_2+\cos\beta\,S_1+\sin\beta\,S_2$ for some $\alpha,\beta\in\mb{R}$ to see
that $X=\cos(\beta-\alpha)\,JX+\sin(\beta-\alpha)\,KX$. The condition \eqref{condition} gives
\begin{equation*}
\tan^2(\beta-\alpha)=\frac{\sin^2(\beta-\alpha)}{\cos^2(\beta-\alpha)}=\left(\frac{\theta_2}{\theta_1}\right)^2
=-\frac{\mu_0+3\mu_1}{\mu_0+3\mu_2}\cdot \frac{\mu_2}{\mu_1},
\end{equation*}
so we can take $\alpha=0$, $\beta=\arctan\sqrt{-\frac{\mu_0+3\mu_1}{\mu_0+3\mu_2}\cdot \frac{\mu_2}{\mu_1}}$ 
and apply Theorem \ref{nontot} to get an anti-Clifford $R=\mu_0 R^0+\mu_1 R^{J}+\mu_2 R^{K}$ which is not totally Jacobi-dual.

\section*{Acknowledgement}

Partially supported by the Serbian Ministry of Education and Science, project 174012.

%\section*{References}

%\bibliography{mybibfile}

\bibliographystyle{amsplain}

\end{document}